\documentclass[11pt]{article}
\usepackage[utf8]{inputenc}
\usepackage[english]{babel}
\usepackage{amssymb,amsmath}
\usepackage[usenames, dvipsnames]{color}
\usepackage{amsthm}
\usepackage{cite}
\let\OLDthebibliography\thebibliography
\renewcommand\thebibliography[1]{
  \OLDthebibliography{#1}
  \setlength{\parskip}{0pt}
  \setlength{\itemsep}{0pt plus 0.3ex}
}
\newtheorem{theorem}{{\scshape Theorem}}[section]
\newtheorem{lemma}[theorem]{{\scshape Lemma	}}
\newtheorem{corollary}[theorem]{{\scshape Corollary}}
\newtheorem{proposition}[theorem]{{\scshape Proposition}}

\newtheorem{example}[theorem]{{\scshape Example}}

\begin{document}
\title{Connections between properties of the additive and the multiplicative groups of a two-sided skew brace}
\author{Timur Nasybullov\footnote{KU Leuven KULAK, 
Etienne Sabbelaan 53, 8500 Kortrijk, Belgium, 
 timur.nasybullov@mail.ru}}
\date{}
\maketitle
\begin{abstract}
We study relations between the additive and the multiplicative groups of a two-sided skew brace. In particular, we prove that if the additive group of a two-sided skew brace is finite solvable (respectively, finitely generated nilpotent, finitely generated residually nilpotent, finitely generated residually finite), then the multiplicative group of this skew brace is solvable (respectively, solvable, residually solvable, residually finite). Also we prove that if the multiplicative group of a two-sided skew brace is nilpotent of nilpotency class $k$, then the additive group of this skew brace is solvable of class at most $2k$. The letter result generalizes the result from \cite[Theorem~A.9]{SmoVen} which says that if the multiplicative group of a finite skew brace  is abelian, then the additive group of this skew brace is solvable.

In addition we solve two problems (Problem 19.49 and Problem 19.90(a)) concerning skew braces which are formulated in the Kourovka notebook.\\

\end{abstract}
\section{Introduction}
\textit{A skew brace} $A=(A,\oplus,\odot)$ is an algebraic system with two binary algebraic operations  $\oplus$, $\odot$ such that $A_{\oplus}=(A,\oplus)$, $A_{\odot}=(A,\odot)$ are groups and the equality
\begin{equation}\label{mainleft}
a\odot(b\oplus c)=(a\odot b)\ominus a \oplus (a\odot c)
\end{equation}
holds for all $a,b,c\in A$, where $\ominus a$ denotes the inverse to $a$ element with respect to the operation~$\oplus$ (we denote by $a^{-1}$ the inverse to $a$ element with respect to $\odot$). The group $A_{\oplus}$ is called \textit{the additive group of a skew brace} $A$, and the group $A_{\odot}$ is called \textit{the multiplicative group of a skew brace} $A$. If $A_{\oplus}$ is abelian, then $A$ is called \textit{a classical brace}.

Classical braces were introduced by Rump in \cite{Rum} in order to study non-degenerate involutive
set-theoretic solutions of the Yang-Baxter equation. Recall that \textit{a set-theoretical solution of the Yang-Baxter equation} is a pair $(X, r)$,
where $X$ is a set and $r:X\times X\to X\times X$ is a bijective map such that 
$$(r\times id)(id \times r)(r\times id)=(id\times r)(r \times id)(id\times r).$$
The solution $(X,r)$ is said to be \textit{non-degenerate} if for $r(x,y)=(\sigma(x,y),\tau(x,y))$ the maps $\sigma(x,\cdot),\tau(\cdot,y):X\to X$ are bijective for fixed $x,y\in X$, and the solution is said to be \textit{involutive} if $r^2=id$. The Yang-Baxter equation first appeared in theoretical physics and statistical
mechanics in the works of Yang \cite{Yan} and Baxter \cite{Bax1,Bax2} and it has led to several
interesting applications in different fields of mathematics. For example, the Yang-Baxter equation appears in topology and algebra since it is connected  with braid groups. The problem of studying set-theoretical solutions of the Yang-Baxter equation was formulated by Drinfel'd in \cite{Dri}. If $A=(A,\oplus,\odot)$ is a classical brace, then the pair $(A,r)$ for 
\begin{equation}\label{ybe}
r(x,y)=\Big(\ominus x\oplus(x\odot y),(\ominus x\oplus(x\odot y))^{-1}\odot x\odot y\Big)
\end{equation}
is a non-degenerate involutive
set-theoretic solution of the Yang-Baxter equation. So, classical braces are useful for the study of the Yang-Baxter equation.

Skew braces were introduced by Guarnieri and Vendramin in \cite{GuaVen} in order to study non-degenerate   
set-theoretic solutions of the Yang-Baxter equation which are not necessarily involutive. For a given skew brace $A$ one can construct non-degenerate   
set-theoretic solutions of the Yang-Baxter equation using formula $(\ref{ybe})$. Skew braces have connections with other algebraic structures such as groups with exact factorizations, Zappa-Sz\'{e}p products,
triply factorized groups and Hopf-Galois extensions  \cite{SmoVen}. They also have applications in knot theory due to connections with biquandles and racks  \cite{Bac,SmoVen}. Some algebraic aspects of skew braces are studied in \cite{CedSmoVen, SmoVen, KonSmoVen}. A big list of problems concerning skew braces is collected in \cite{Ven}.

Equality (\ref{mainleft}) makes the additive group $A_{\oplus}$ and the multiplicative group $A_{\odot}$ of a skew brace $A$ strongly connected with each other. For example, the following problems formulated in the Kourovka notebook \cite[Problems 19.49, 19.90]{kt} tell about some of such connections.
\begin{enumerate}
\item Let $A$ be a skew brace with left-orderable multiplicative group. Is the additive group of $A$ left-orderable?
\item  Does there exist a skew brace with solvable additive group but non-solvable multiplicative group?
\item  Does there exist a skew brace with nilpotent multiplicative group but non-solvable additive group?
\end{enumerate}
Questions 2 and 3 from the list above are known to have the negative answer in some particular cases. For example, the following result is proved in \cite[Corollary~1.23]{SmoVen}.
\begin{theorem}\label{finnil}Let $A$ be a finite skew brace. If $A_{\oplus}$ is nilpotent, then $A_{\odot}$ is solvable.
\end{theorem}
Another result is formulated in \cite[Theorem A.9]{SmoVen} and proved in \cite[Theorem 2]{Byo} in terms of Hopf-Galois structures. 
\begin{theorem}\label{multab} Let $A$ be a finite skew brace. If $A_{\odot}$ is abelian, then $A_{\oplus}$ is solvable. 
\end{theorem}
In the present text we are going to study the formulated problems for two-sided skew braces.  The skew brace is called \textit{two-sided} if together with equality (\ref{mainleft}) the equality 
\begin{equation}\label{mainright}
(a\oplus b)\odot c=(a\odot c)\ominus c \oplus (b \odot c)
\end{equation}
holds for all $a,b,c\in A$. If the multiplicative group of a skew brace $A$ is abelian, then this skew brace is two-sided. However, not every two-sided skew brace has abelian multiplicative group (see \cite[Example 1.18]{SmoVen}). Two-sided classical braces are in one-to-one correspondence with radical rings \cite{Rum2}. So, the notion of a two-sided skew brace generalizes the notion of a radical ring. In the present paper we prove that if the additive group of a two-sided skew brace is finite solvable (respectively, finitely generated nilpotent, finitely generated residually nilpotent, finitely generated residually finite), then the multiplicative group of this skew brace is solvable (respectively, solvable, residually solvable, residually finite). Also we prove that if the multiplicative group of a two-sided skew brace is nilpotent of nilpotency class $k$, then the additive group of this skew brace is solvable of class at most $2k$. The letter result generalizes Theorem \ref{multab}.

\section{Definitions and known results}
In this section we recall some notions and known facts about skew braces.  Equality (\ref{mainleft}) implies that the unit element of $A_{\oplus}$ coincides with the unit element of $A_{\odot}$, we denote this element by $1$. Also from equality (\ref{mainleft}) follows that the equality 
\begin{equation}\label{1inv}
a\odot(\ominus b\oplus c)=a\ominus(a\odot b)\oplus (a\odot c)
\end{equation}
holds for all $a,b,c\in A$.

For $a\in A$ the map $\lambda_a: x\mapsto \ominus a\oplus (a\odot x)$ is an automorphism of the group $A_{\oplus}$. Moreover, the map $a\mapsto \lambda_a$ gives a homomorphism $A_{\odot}\to {\rm Aut}(A_{\oplus})$. The subset $I$ of $A$ is called \textit{an ideal of} $A$ if $I$ is a normal subgroup of both $A_{\oplus}$, $A_{\odot}$ and $\lambda_a(I)=I$ for all $a\in A$. The letter condition means that $a\oplus I=a\odot I$ for all $a\in A$. It is clear that $I$ is itself a skew brace, so, we can speak about $I_{\oplus}$, $I_{\odot}$. For an ideal $I$ of a skew brace $A$ \textit{the quotient} $A/I$ is a natural skew brace defined on the set of cosets $a\oplus I$ with the operations given by $(a\oplus I)\oplus (b\oplus I)= (a\oplus b)\oplus I$, $(a\oplus I)\odot (b\oplus I)= (a\odot b)\oplus I$. The skew brace $A/I$ has the additive group $(A/I)_{\oplus}=A_{\oplus}/I_{\oplus}$ and multiplicative group $(A/I)_{\odot}=A_{\odot}/I_{\odot}$. A skew brace is called \textit{simple} if it has no proper ideals.

According to \cite{CedSmoVen} for elements $a,b\in A$ denote by $a*b=\ominus a\oplus (a\odot b)\ominus b=\lambda_a(b)\ominus b$. Direct calculations imply the equality 
\begin{equation}\label{distrr}
a*(b+c)=(a*b)\oplus b\oplus (a*c)\ominus b
\end{equation}
for $a,b,c\in A$. For subsets $X, Y\subseteq A$ denote by $X*Y$ the subgroup $X*Y=\langle x*y~|~x\in X, y\in Y \rangle_{\oplus}$ of $A_{\oplus}$ additively generated by elements $x*y$ for $x\in X, y\in Y$. The following result is proved in \cite[Proposition~2.1]{CedSmoVen}
\begin{proposition}\label{a*a} Let $A$ be a skew brace. Denote by $A^{(1)}=A$, $A^{(k+1)}=A^{(k)}*A$ for $k>1$. Then $A^{(k)}$ is an ideal of $A$ for all $k$.
\end{proposition}
A skew brace $A$ is called \textit{trivial} if $a\oplus b=a\odot b$ for all $a,b\in A$. The following statement is proved in \cite[Proposition 2.3]{CedSmoVen}.
\begin{proposition}\label{fac}Let $A$ be a skew brace. Then $A^{(2)}$ is the smallest ideal of $A$ such that $A/A^{(2)}$ is a trivial skew brace.
\end{proposition}
\section{Examples and counterexamples}\label{examples}
In this section we give examples of skew braces which give answers to several problems.
\begin{example}{\rm On the set of integers consider two operations $\oplus$ and $\odot$ given by the formulas.
\begin{align}
\notag  a\oplus b=a+(-1)^ab, && a\odot b=a+b
\end{align}
Then $A=(\mathbb{Z},\oplus,\odot)$ is a skew brace with infinite cyclic multiplicative group $A_{\odot}=\mathbb{Z}$ and infinite dihedral additive group $A_{\oplus}=\mathbb{Z}\rtimes_{-id}\mathbb{Z}_2$.}
\end{example}
 This skew brace gives negative answers to the following two problems.
\medskip

\noindent\textbf{\cite[Question A.10]{SmoVen}.} Let $A$ be a skew brace with multiplicative group isomorphic to $(\mathbb{Z},+)$. Is it trues that the additive group of $A$ is also isomorphic to $(\mathbb{Z},+)$?
\medskip

The group $A_{\oplus}$ is not abelian, therefore the answer to \cite[Question A.10]{SmoVen} is negative.
\medskip

\noindent\textbf{\cite[Problem 19.49]{kt}.} Let $A$ be a skew brace with left-orderable multiplicative group. Is the additive
group of $A$ left-orderable?
\medskip

The group $A_{\odot}$ is infinite cyclic and therefore is left-orderable group. At the same time the group $A_{\oplus}$ has $2$-torsion, therefore $A_{\oplus}$ is not left-orderable and the answer to \cite[Problem 19.49]{kt} is negative.

Recently the preprint \cite{CedSmoVen} appeared, where the same example is constructed \cite[Theorem~5.8]{CedSmoVen}. In this preprint it is noticed that this example  gives the answer to \cite[Question A.10]{SmoVen}. Since this example is not too big and we constructed it before reading \cite{CedSmoVen}, we decided not to delete it from the text and notice that it  also solves \cite[Problem 19.49]{kt}.
\begin{example}\label{exxx}{\rm For $n\geq 2$ denote by $U_n$ the set of strictly upper triangular matrices of degree $n$ over $\mathbb{Z}$. For $A,B\in U_n$ denote by $A\oplus B$, $A\odot B$ the matrices from $U_n$ of the following form
\begin{align}
\notag A\oplus B=A+B,&& A\odot B=({\rm I}_n+A)({\rm I}_n+B)-{\rm I}_n,
\end{align}
where ${\rm I}_n$ is the identity matrix. It is obvious that $(U_n,\oplus)$ is the group isomorphic to $\mathbb{Z}^{n(n-1)/2}$, and $(U_n,\odot)$ is the group isomorphic to the group ${\rm UT}_n(\mathbb{Z})$ of upper unitriangular matrices over $\mathbb{Z}$. For arbitrary matrices $A,B,C\in U_n$ we have
\begin{align}
\notag A\odot(B\oplus C)&=A\odot(B+C)\\
\notag&=({\rm I}_n+A)({\rm I}_n+B+C)-{\rm I_n}\\
\notag&=({\rm I}_n+A)({\rm I}_n+B)-{\rm I_n}+ ({\rm I}_n+A)C\\
\notag&=({\rm I}_n+A)({\rm I}_n+B)-{\rm I_n}+ ({\rm I}_n+A)(I_n+C)-(I_n+A)\\
\notag&=\big(({\rm I}_n+A)({\rm I}_n+B)-{\rm I_n}\big)-A+\big(({\rm I}_n+A)(I_n+C)-I_n\big)\\
\notag&=(A\odot B)\ominus A\oplus (A\odot C),
\end{align}
therefore $A_n=(U_n,\oplus,\odot)$ is a skew brace with $(A_n)_{\oplus}=\mathbb{Z}^{n(n-1)/2}$, $(A_n)_{\odot}={\rm UT}_n(\mathbb{Z})$.}
\end{example}

Let $A$ be a skew brace which is a direct sum of skew braces $A_n$ from Example \ref{exxx} for $n=2,3,\dots$, i.~e. $A=\{a=(a_2,a_3,\dots)~|~a_i\in A_i, |{\rm supp}(a)|<\infty\}$ and operations are componentwise. The additive group $A_{\oplus}$ is isomorphic to the direct sum of infinite number of copies of $\mathbb{Z}$, i.~e. $A_{\oplus}$ is abelian. The multiplicative group $A_{\odot}$ is isomorphic to the direct sum of groups ${\rm UT}_n(\mathbb{Z})$ for $n=2,3,\dots$ Since the group ${\rm UT}_n(\mathbb{Z})$ is solvable of degree $\lceil{\rm log}_2(n)\rceil$, the group $A_{\odot}$ is not solvable. The skew brace $A$ gives an answer to the following problem.
\medskip

\noindent\textbf{\cite[Problem 19.90(a)]{kt}.} Does there exist a skew brace with solvable additive group but non-solvable multiplicative group?
\begin{example}\label{exx3}{\rm For $n\geq2$ denote by $T_n$ the set of invertible upper triangular matrices of degree $n$ over $\mathbb{Z}$. Since every matrix $X$ from $T_n$ can be uniquely expressed as a  product $X=({\rm I}_n+A)a$ of the upper unitriangular matrix $({\rm I}_n+A)$ and diagonal matrix $a$, the set $T_n$ can be considered as the set of pairs $T_n=\{(A,a)~|~A\in U_n, a\in D_n\}$, where $D_n$ is the set of diagonal matrices of degree $n$ with $\pm 1$ on the diagonal. For $(A,a), (B,b)\in T_n$ denote by $(A,a)\oplus (B,b)$, $(A,a)\odot (B,b)$ the following matrices from $T_n$
\begin{align}
\notag(A,a)\oplus (B,b)=(A+B, ab),&&(A,a)\odot (B,b)=\Big((A+{\rm I_n})a(B+I_n)a^{-1}-I_n, ab\Big).
\end{align}
Using direct calculations it is easy to see that $A_n=(T_n,\oplus,\odot)$ is a skew brace. The additive group $A_{\oplus}$ is isomorphic to $\mathbb{Z}^{n(n-1)/2}\oplus\mathbb{Z}_2^n$, and the multiplicative group  $A_{\odot}$ is isomorphic to the group ${\rm T}_n(\mathbb{Z})$ of upper triangular matrices of degree $n$ over $\mathbb{Z}$ with isomorphism given by $(A,a)\mapsto ({\rm I}_n+A)a$.}
\end{example}
The direct sum of skew braces $A_n$ from Example \ref{exx3} for $n=2,3,\dots$ gives another example of a skew brace which answers \cite[Problem 19.90(a)]{kt}.

\section{Two-sided skew braces}
In this section we study connections between properties of the additive and the multiplicative groups of two-sided skew braces. The following simple technical lemma is true for two-sided skew braces.
\begin{lemma}\label{multaut} If $A$ is a two sided skew brace, then the equality
$$c^{-1}\odot(a\oplus b)\odot c=(c^{-1}\odot a\odot c)\oplus (c^{-1}\odot b\odot c)$$
holds for all $a,b,c\in A$. 
\end{lemma}
\noindent \textbf{Proof.} Applying equalities (\ref{mainleft}), (\ref{1inv}) and (\ref{mainright}) using direct calculations we have
\begin{align}
\notag c^{-1}\odot(a\oplus b)\odot c&\overset{(\ref{mainright})}{=}c^{-1}\odot\big((a\odot c)\ominus c\oplus (b\odot c)\big)\\
\notag &\overset{(\ref{mainleft})}{=}(c^{-1}\odot a\odot c)\ominus c^{-1} \oplus c^{-1}\big(\ominus c\oplus (b\odot c)\big)\\
\notag &\overset{(\ref{1inv})}{=}(c^{-1}\odot a\odot c)\ominus c^{-1} \oplus c^{-1}\ominus (c^{-1}\odot c) \oplus (c^{-1}\odot b\odot c)\\
\notag &=(c^{-1}\odot a\odot c) \oplus (c^{-1}\odot b\odot c),
\end{align}
the lemma is proved.\hfill $\square$

Lemma \ref{multaut} has the following corollary which answers the question \cite[Problem 38]{Ven} in the case when $A$ is a finite skew brace.
\begin{corollary}Let $A$ be a finite two-sided skew brace with solvable additive group. Then $A$ is simple if and only if $A$ is a trivial skew brace of prime order $p$.
\end{corollary}
\noindent \textbf{Proof.} If $A$ is a trivial skew brace of prime order $p$, then it is obviously simple, and we need to prove only the if part of the corollary. 
From Lemma \ref{multaut} follows that the map $a\mapsto c^{-1}\odot a\odot c$  is an automorphism of $A$ (and therefore an automorphism of $A_{\oplus}$), therefore every characteristic subgroup of $A_{\oplus}$ is an ideal of $A$. Denote by $I$ the commutator subgroup of $A_{\oplus}$. Since $A_{\oplus}$ is solvable, $I\neq A$. Since $I$ is characteristic in $A_{\oplus}$, it is an ideal of $A$, and since $A$ is simple, $I=1$ and $A_{\oplus}$ is abelian. For a prime $p$ dividing $|A|$ denote by $J$ the Sylow $p$-subgroup of $A_{\oplus}$. Since $A_{\oplus}$ is abelian, there is only one Sylow $p$-subgroup in $A_{\oplus}$, therefore $J$ is characteristic in $A_{\oplus}$ and $J$ is an ideal of $A$. Since $A$ is simple, $A=J$ and $|A|=p^n$ for some $n$, so, the corollary follows from \cite[Corollary 4.8]{CedSmoVen}.\hfill$\square$

If $A$ is a two sided skew brace, then Theorem \ref{finnil} can be generalized in the following form. 
\begin{theorem} Let $A$ be a two-sided skew brace. Then
\begin{enumerate}
\item If $A$ is finite and $A_{\oplus}$ is solvable, then $A_{\odot}$ is solvable,
\item If $A_{\oplus}$ is finitely generated and nilpotent, then $A_{\odot}$ is solvable,
\item If $A_{\oplus}$ is finitely generated and residually nilpotent, then $A_{\odot}$ is residually solvable,    
\item If $A_{\oplus}$ is finitely generated and residually finite, then $A_{\odot}$ is residually finite.
\end{enumerate}
\end{theorem}
\noindent \textbf{Proof.} 1) If $A_{\oplus}$ is abelian, then the result follows from Theorem \ref{finnil}. If $A_{\oplus}$ is not abelian, then denote by $I$ the commutator subgroup of $A_{\oplus}$. Since $I\neq1$ is a characteristic subgroup of $A_{\oplus}$, it is an ideal of $A$. The quotient $A/I$ has a solvable additive group $(A/I)_{\oplus}$. Therefore using induction on the order of a skew brace we can assume that $(A/I)_{\odot}$ is solvable. Also since $I_{\oplus}$ is solvable, using induction on the order of a skew brace we can assume that $I_{\odot}$ is solvable. Therefore $A_{\odot}$ is solvable.

2) If $A_{\oplus}$ is finitely generated abelian, then $A_{\oplus}$ is a direct product of a finite abelian group $B$ and a free abelian group $\mathbb{Z}^n$ for some $n$. Since $B$ is a characteristic subgroup of $A_{\oplus}$, it is an ideal of $A$ (in particular, $B$ is a skew brace). Since $B_{\oplus}$ is finite and abelian, from Theorem \ref{finnil} follows that $B_{\odot}$ is solvable.

Let $p>|B|$ be a prime number. Denote by $I_p$ the subgroup of $A_{\oplus}$ generated by $p$-th powers (with respect to $\oplus$) of all elements. Since $I_p$ is a characteristic subgroup of $A_{\oplus}$, it is an ideal of $A$. The quotient $A/I_p$ has the order $p^n$. Therefore the multiplicative group $(A/I_p)_{\odot}$ is nilpotent of nilpotency class at most $n$ and for arbitrary elements $x_1, \dots, x_{n+1}\in A$ we have $[x_1,\dots,x_{n+1}]_{\odot}\oplus I_p=I_p$, where $[x_1,\dots,x_{n+1}]_{\odot}$ denotes the simple commutator of length $n+1$ of elements $x_1,\dots,x_{n+1}$ with respect to $\odot$. Therefore for all elements $x_1,\dots,x_{n+1}$ the commutator $[x_1,\dots,x_{n+1}]_{\odot}$ belongs to $I_p$ for all $p>|B|$. Hence the $n$-th derived subgroup $A_{\odot}^{(n)}$ of $A_{\odot}$ belongs to $\cap_{p>|B|} I_p=B$ (this equality follows from the fact that $p>|B|$). Since $B_{\odot}$ is solvable, $A_{\odot}$ is solvable.

If $A_{\oplus}$ is nilpotent of class $k>1$, then denote by $I$ the smallest nontrivial term of the lower central series of $A_{\oplus}$. Since $I$ is a characteristic subgroup of $A_{\oplus}$, it is an ideal of $A$. The quotient $A/I$ has a finitely generated nilpotent additive group $(A/I)_{\oplus}$ of nilpotency class $k-1$. Therefore using induction on the nilpotency class of $A_{\oplus}$ we can assume that $(A/I)_{\odot}$ is solvable. Since $A_{\oplus}$  is finitely generated nilpotent, $I_{\oplus}$ is finitely generated abelian, and by the paragraph above $I_{\odot}$ is solvable. Therefore $A_{\odot}$ is solvable.

3) Let $x\in A$. Denote by $I$ the minimal member of the lower central series of $A_{\oplus}$ which does not contain $x$ (such member exists since $A_{\oplus}$ is residually nilpotent). Since $I$ is characteristic subgroup of $A_{\oplus}$, it is an ideal of $A$. The quotient $A/I$ has finitely generated nilpotent additive group $(A/I)_{\oplus}$, therefore from 2) follows that $(A/I)_{\odot}$ is solvable.

4) Since $A_{\oplus}$ is residually finite, for an element $x\neq 1$ there exists a normal subgroup $N\triangleleft A_{\oplus}$ of index $n$ such that $x\notin N$. Since $A_{\oplus}$ is finitely generated, there exists only a finite number of subgroups of index $n$, therefore there is a characteristic subgroup $I=\cap_{\varphi\in{\rm Aut}(A_{\oplus})} \varphi(N)$ of $A_{\oplus}$ of finite index such that $x\notin I$. Since $I$ is characteristic, it is an ideal of $A$, therefore it is a normal subgroup of $A_{\odot}$.\hfill $\square$

The following lemma can be thought as a stronger version of Proposition \ref{a*a} which is true for two-sided skew braces. 
\begin{lemma}\label{2s2} Let $A$ be a two-sided skew brace and $X$ be a normal subgroup of $A_{\odot}$. Then $X*A$ is an ideal of $A$.
\end{lemma}
\noindent \textbf{Proof.} For arbitrary elements $x\in X$, $a,b\in A$ we have the equality
$$
 b\oplus (x*a)\ominus b=b\oplus \lambda_x(a)\ominus a\ominus b=b\ominus \lambda_x(b)\oplus\lambda_x(b\oplus a)\ominus(b\oplus a)=\ominus(x*b)\oplus(x*(b\oplus a)),$$
therefore $(X*A)_{\oplus}$ is a normal subgroup of $A_{\oplus}$. Again for arbitrary elements $x\in X$, $a,b\in A$ we have
\begin{align}
\label{1conj}\lambda_b(x*a)&=\lambda_b(\lambda_x(a)\ominus a)=\lambda_b\lambda_x(a)\ominus \lambda_b(a)=\lambda_b\lambda_x\lambda_b^{-1}(\lambda_b(a))\ominus \lambda_b(a)
=(b\odot x\odot b^{-1})*\lambda_b(a),
\end{align}
therefore $\lambda_b(X*A)=X*A$ for all elements $b\in A$. For elements $x,y \in X$, $a,b\in A$ we have $(x*a)\odot(y*b)=(x*a)\oplus\lambda_{x*a}(y*b)$, therefore $X*A$ is closed under multiplication and we can speak about $(X*A)_{\odot}$. Finally, using Lemma \ref{multaut} for arbitrary elements $x\in X$, $a,b\in A$ we have
\begin{align}
\notag b^{-1}\odot(x*a)\odot b&=b^{-1}\odot(\ominus x\oplus(x\odot a)\ominus a)\odot b\\
\notag&=\ominus (b^{-1}\odot x\odot b)\oplus (b^{-1}\odot x\odot a\odot b)\ominus (b^{-1}\odot a\odot b)\\
\notag&=(b^{-1} \odot x\odot b)*(b^{-1} \odot a\odot b),
\end{align}
therefore $(X*A)_{\odot}$ is a normal subgroup of $A_{\odot}$ and $X*A$ is an ideal of $A$.\hfill$\square$
\begin{lemma}\label{2s2} Let $A$ be a two-sided skew brace and $Z$ be the center of $A_{\odot}$. Then $I=(A*Z)\oplus (Z*A)$ is an ideal of $A$ such that the group $I_{\oplus}$ is abelian.
\end{lemma}
\noindent \textbf{Proof.} We need to check that $I$ is a normal subgroup of both $A_{\oplus}$, $A_{\odot}$ and that $\lambda_a(I)=I$ for all $a\in A$. By Lemma \ref{2s2}, the set $Z*A$ is an ideal of $A$. Using this fact, for $a\in A$, $x\in Z$ we have
\begin{align}
\notag a*x&=\ominus a\oplus (a\odot x)\ominus x=\ominus a\oplus x\ominus x\oplus (x\odot a)\ominus a\oplus a\ominus x\\
\label{modd}&=\ominus a\oplus x\oplus(x*a)\oplus a\ominus x= \ominus a\oplus x\oplus a\ominus x\oplus c
\end{align}
for $c=(x\ominus a)\oplus (x*a)\ominus(x\ominus a)\in Z*A$. For arbitrary elements $a,b\in A$, $x\in Z$ using equality (\ref{modd}) several times we have
\begin{align}
\notag\ominus b\oplus (a*x)\oplus b&\overset{(\ref{modd})}{=} \ominus b\ominus a\oplus x\oplus a\ominus x\oplus c\oplus b&c\in& Z*A,\\
\notag&= \ominus b\ominus a\oplus x\oplus a\ominus x\oplus b\oplus c_1&c_1\in& Z*A,\\
\notag&= \ominus (a\oplus b)\oplus x\oplus (a\oplus b)\ominus x\oplus x\ominus b\ominus x\oplus b\oplus c_1&c_1\in& Z*A,\\
\notag&\overset{(\ref{modd})}{=} ((a\oplus b)*x)\oplus c_2\ominus ((b*x)\oplus c_3)\oplus c_1&c_1,c_2,c_3\in& Z*A,\\
\notag&= ((a\oplus b)*x)\ominus (b*x)\oplus c_4&c_4\in& Z*A,
\end{align}
therefore $\ominus b\oplus (A*Z)\oplus b\subseteq I$ for arbitrary $b\in A$, and $I_{\oplus}$ is a normal subgroup of $A_{\oplus}$. For arbitrary elements $a,b\in A$, $x\in Z$ we have
\begin{align}
\notag(b\odot a\odot b^{-1})\odot(\ominus b\oplus x\oplus b)&\overset{(\ref{1inv})}{=}(b\odot a\odot b^{-1})\ominus (b\odot a)\oplus ((b\odot a\odot b^{-1})\odot (x\oplus b))\\
\notag&\overset{(\ref{mainleft})}{=}(b\odot a\odot b^{-1})\ominus (b\odot a)\oplus (b\odot a\odot b^{-1}\odot x)\ominus (b\odot a\odot b^{-1})\\
\notag &~~~~\oplus (b\odot a)\\
\notag&=(b\odot a\odot b^{-1})\ominus (b\odot a)\oplus x\oplus (x*(b\odot a\odot b^{-1}))\oplus (b\odot a)\\
\notag &=(b\odot a\odot b^{-1})\ominus (b\odot a)\oplus x\oplus (b\odot a)\oplus c
\end{align}
for $c=\ominus (b\odot a)\oplus (x*(b\odot a\odot b^{-1}))\oplus (b\odot a)\in Z*A$. From this equality and equality (\ref{modd}) follows that
\begin{align}
\notag(b\odot a\odot b^{-1})*(\ominus b\oplus x\oplus b)&=\ominus (b\odot a)\oplus x\oplus (b\odot a)\ominus (\ominus b\oplus x\oplus b)\oplus c_1&c_1\in& Z*A\\
\notag&=\ominus (b\odot a)\oplus x\oplus (b\odot a) \ominus b\ominus x\oplus b\oplus c_1&c_1\in& Z*A\\
\notag&\overset{(\ref{modd})}{=}((b\odot a)*x)\oplus c_2\ominus (b*x)\oplus c_3\oplus c_1&c_1,c_2,c_3\in& Z*A\\
\notag&=((b\odot a)*x)\ominus (b*x)\oplus c&c\in& Z*A	
\end{align}
Therefore for $x\in Z$, $a,b\in A$ the element $(b\odot a\odot b^{-1})*(\ominus b\oplus x\oplus b)$ belongs to $I$. For arbitrary elements $b\in A$, $x\in Z$ we have $b\oplus\lambda_b(x)=b\odot x=x\odot b=x\oplus \lambda_x(b)=x\oplus (x*b)\oplus b$, therefore 
\begin{equation}\label{lam}
\lambda_b(x)=\ominus b\oplus x\oplus b\oplus c
\end{equation}
for $c=\ominus b\oplus (x*b)\oplus b\in Z*A$. If in equality (\ref{1conj}) we change $x$ and $a$, then we have the equality 
\begin{equation}\label{lambdda}\lambda_b(a*x)=(b\odot a\odot b^{-1})*\lambda_b(x).
\end{equation} 
So, summarizing together equalities (\ref{distrr}), (\ref{lam}), (\ref{lambdda}) and the fact that $(b\odot a\odot b^{-1})*(\ominus b\oplus x\oplus b)$ belongs to $I$ for arbitrary elements $a,b\in A$, $x\in Z$ we have
\begin{align}
\notag\lambda_b(a*x)&\overset{(\ref{lambdda})}{=}(b\odot a\odot b^{-1})*\lambda_b(x)&&\\
\notag &\overset{(\ref{lam})}{=}(b\odot a\odot b^{-1})*(\ominus b\oplus x\oplus b\oplus c)&c\in& Z*A\\
\notag &\overset{(\ref{distrr})}{=}(b\odot a\odot b^{-1})*(\ominus b\oplus x\oplus b)&&\\
\notag &~~~~\oplus(\ominus b\oplus x\oplus b)\oplus ((b\odot a\odot b^{-1})*c)\ominus (\ominus b\oplus x\oplus b)&c\in& Z*A\\
\notag &=(b\odot a\odot b^{-1})*(\ominus b\oplus x\oplus b)&&\\
\notag &~~~~\ominus b\oplus x\oplus b\ominus(b\odot a\odot b^{-1}) \oplus ((b\odot a\odot b^{-1})\odot c)\ominus c\ominus b\ominus x\oplus b&c\in& Z*A\\
\notag &=(b\odot a\odot b^{-1})*(\ominus b\oplus x\oplus b)&&\\
\notag &~~~~\ominus b\oplus x\oplus b\ominus(b\odot a\odot b^{-1}) \oplus (b\odot a\odot b^{-1})\oplus c_1\ominus c\ominus b\ominus x\oplus b&c,c_1\in& Z*A\\
\notag &=(b\odot a\odot b^{-1})*(\ominus b\oplus x\oplus b)\oplus c_2&c_2\in& Z*A,
\end{align}
therefore for arbitrary element $b\in A$ we have $\lambda_b(A*Z)\subseteq I$ and then $\lambda_b(I)\subseteq I$. Finally, from Lemma~\ref{multaut} follows that for arbitrary elements $x\in X$, $a,b\in A$ we have
\begin{align}
\notag b^{-1}\odot(a*x)\odot b&=(b^{-1} \odot a\odot b)*(b^{-1} \odot x\odot b),
\end{align}
therefore $(A*Z)_{\odot}$ is a normal subgroup of $A_{\odot}$, $I_{\odot}$ is a normal subgroup of $A_{\odot}$ and therefore $I$ is an ideal of $A$.

For arbitrary elements  $a,b,c,d\in A$, using equalities (\ref{mainleft}) and  (\ref{mainright}) we have the following equality
\begin{align}
\notag(a\oplus b)\odot(c\oplus d)&\overset{(\ref{mainleft})}{=}((a\oplus b)\odot c)\ominus (a\oplus b)\oplus ((a\oplus b)\odot d)\\
\notag&\overset{(\ref{mainright})}{=}(a\odot c)\ominus c\oplus (b\odot c)\ominus b\ominus a\oplus (a\odot d)\ominus d\oplus (b\odot d)\\
\label{eq4q1}&=(a\odot c)\ominus c\oplus (b\odot c)\ominus b\oplus(a*d)\oplus (b\odot d).
\end{align}
From the other hand, applying equalities  (\ref{mainleft}) and  (\ref{mainright}) in the another sequence we have the following equality 
\begin{align}
\notag(a\oplus b)\odot(c\oplus d)&\overset{(\ref{mainright})}{=}(a\odot(c\oplus d))\ominus (c\oplus d)\oplus (b\odot(c\oplus d))\\
\notag&\overset{(\ref{mainleft})}{=}(a\odot c)\ominus a\oplus (a\odot d)\ominus d\ominus c\oplus (b\odot c)\ominus b\oplus (b\odot d)\\
\label{eq5q1}&=(a\odot c)\oplus (a*d)\ominus c\oplus (b\odot c)\ominus b\oplus (b\odot d).
\end{align}
Since expression (\ref{eq4q1}) and (\ref{eq5q1}) represent the same element of $A$, they must be equal, this fact means that  for all $a,b,c,d\in A$ we have the equality
\begin{equation}\label{nahodka}
(\ominus c\oplus (b\odot c)\ominus b)\oplus(a*d)=(a*d)\oplus (\ominus c\oplus (b\odot c)\ominus b).
\end{equation}
Applying equality (\ref{nahodka}) to elements $a,c\in Z$, $b,d\in A$, we conclude that $(Z*A)_{\oplus}$ is abelian; applying this equality to elements $b,d\in Z$, $a,c\in A$, we conclude that $(A*Z)_{\oplus}$ is abelian; and applying the same equality to elements $a,b\in Z$, $c,d\in A$ we conclude that every element from  $(A*Z)_{\oplus}$ commutes (with respect to $\oplus$) with every element from $(Z*A)_{\oplus}$, therefore $I_{\oplus}$ is abelian group.
\hfill$\square$

\begin{theorem}\label{mgener} Let $A$ be a two-sided skew brace. If $A_{\odot}$ is nilpotent of nilpotency class $k$, then $A_{\oplus}$ is solvable of class at most $2k$. 
\end{theorem}
\noindent \textbf{Proof.} We will prove the statement using induction on the nilpotency class $k$ of $A_{\odot}$. If $k=1$, then the group $A_{\odot}$ is abelian. By Proposition \ref{fac} the skew brace $A/(A*A)$ is trivial, i.~e. $a\odot b=a\oplus b$ for all $a,b \in A/(A*A)$, therefore $(A/(A*A))_{\oplus}$ is abelian. By Lemma \ref{2s2} the group $(A*A)_{\oplus}$ is abelian, so, $A_{\oplus}$ has a normal abelian subgroup $(A*A)_{\oplus}$ with abelian quotient $A_{\oplus}/(A*A)_{\oplus}$, i.~e. $A_{\oplus}$ is solvable of class at most $2$. The basis of induction ($k=1$) is proved.

Let $k>1$, denote by $Z$ the center of $A_{\odot}$ and let $I=(A*Z)\oplus (Z*A)$. By Lemma \ref{2s2} the set $I$ is an ideal of $A$ with abelian group $I_{\oplus}$.
For all $x\in Z$, $a\in A$ there exist some elements $c,d\in I$ such that $x\odot a=x\oplus a\oplus c$, $a\odot x=a\oplus x\oplus d$, therefore the set $J=Z\oplus I$ is the ideal of $A$. Since $Z\subseteq J$, the multiplicative group $(A/J)_{\odot}$ of $A/J$ is nilpotent of nilpotency class at most $k-1$, therefore using induction hypothesis we can assume that $(A/J)_{\oplus}$ is solvable of class at most $2(k-1)$. So, we have a normal series of groups $A_{\oplus}>J_{\oplus}>I_{\oplus}>1$ such that $A_{\oplus}/J_{\oplus}$ is solvable of class at most $2k-2$, $J_{\oplus}/I_{\oplus}$ and $I_{\oplus}$ are abelian, therefore $A_{\oplus}$ is solvable of class at most $2k$.\hfill$\square$

If the multiplicative group of a skew brace is abelian, then this skew brace is two-sided. The following result follows from Theorem \ref{mgener} and generalizes Theorem \ref{multab}.
\begin{corollary}\label{nill}Let $A$ be a skew brace. If  $A_{\odot}$ is abelian, then $A_{\oplus}$ is metabelian.
\end{corollary}

Corollary \ref{nill} generalizes Theorem \ref{multab} in two directions: at first, Corollary \ref{nill} is true for skew braces of arbitrary cardinality (not only for finite skew braces),  at second, it gives the upper bound of the solvability degree of $A_{\oplus}$ ($A_{\oplus}$ must be metabelian). Also the proof of Theorem \ref{mgener} (which implies Corollary \ref{nill}) is elementary, i.~e. it does not use classification of finite simple groups and other advanced results from finite group theory which are used in the proof of \cite[Theorem 2]{Byo} which implies Theorem \ref{multab}.


{\small
\begin{thebibliography}{00}
\bibitem{Bac}
D.~Bachiller, Solutions of the Yang-Baxter equation associated to skew left braces, J.~Knot Theory Ramifications V.~27, N.~8, 2018, 1850055.
\bibitem{Bax1}
 R.~Baxter, Partition function of the eight-vertex lattice model, Ann. Physics, V.~70, 1972, 193--228.
\bibitem{Bax2}
R.~Baxter, Exactly solved models in statistical mechanics, Academic Press, Inc., London, 1982.
\bibitem{Byo}
N.~Byott, Solubility criteria for Hopf-Galois structures, New York J. Math., V.~21, 2015, 883--903.
\bibitem{CedSmoVen}
F.~Ced\'o, A.~Smoktunowicz, L.~Vendramin, Skew left braces of nilpotent type, 	ArXiv:Math/1806.01127.
\bibitem{Dri}
 V.~Drinfel'd, On some unsolved problems in quantum group theory, Quantum groups (Leningrad, 1990), 1-8, 
Lecture Notes in Math., 1510, Springer, Berlin, 1992. 
\bibitem{GuaVen}
L. Guarnieri, L. Vendramin, Skew braces and the Yang-Baxter equation, Math. Comp., V.~86, N.~307, 2017, 2519--2534.
\bibitem{KonSmoVen}
A.~Konovalov, A.~Smoktunowicz, L.~Vendramin, On skew braces and their ideals, ArXiv:Math/1804.04106.
\bibitem{kt}
The Kourovka notebook, Unsolved problems in group theory. Edited by V.~D.~Mazurov and E.~I.~Khukhro, 19-th. ed.. Russian Academy of Sciences Siberian Division. Institute of Mathematics, Novosibirsk, 2018.
\bibitem{Rum}
W.~Rump, Braces, radical rings, and the quantum Yang-Baxter equation, J.~Algebra, V.~307, N.~1, 2007, 153--170.
\bibitem{Rum2}
 W.~Rump, Modules over braces, Algebra Discrete Math., N.~2, 2006, 127--137.
\bibitem{SmoVen}
A.~Smoktunowicz, L.~Vendramin, On skew braces (with an appendix by N.~Byott and L.~Vendramin), Journal of combinatorial algebra, V.~2, N.~1, 47--86.
\bibitem {Ven}
L.~Vendramin, Problems on skew left braces,  ArXiv:Math/1807.06411.
\bibitem{Yan}
C.~Yang, Some exact results for the many-body problem in one dimension with repulsive
delta-function interaction, Phys. Rev. Lett., V.~19, 1967, 1312--1315.
\end{thebibliography}
}
\end{document}